\documentstyle[twoside,12pt]{article}
\newtheorem{prop}{Proposition}[section]
\newtheorem{dfn}[prop]{Definition}
\newtheorem{theo}[prop]{Theorem}

\newtheorem{rem}[prop]{Remark}
\newtheorem{coro}[prop]{Corollary}

\newtheorem{exam}[prop]{Example}

\title{On Calabi-Yau Complete Intersections
in  Toric Varieties}

\author{Victor V. Batyrev\thanks{Supported by DFG.} \\
\small  Universit\"at-GHS-Essen, Fachbereich  6,  Mathematik \\
\small  45141  Essen, Germany  \\
\small e-mail: victor.batyrev@aixrs1.hrz.uni-essen.de \\
and \\
Lev A. Borisov \\
\small Department of Mathematics,  University of Michigan \\
\small Ann Arbor, Michigan 48109-1003, USA \\
\small e-mail: lborisov@math.lsa.umich.edu}

\begin{document}

\date{}

\maketitle

\begin{abstract}
We investigate  Hodge-theoretic  properties
of Calabi-Yau complete intersections $V$ of $r$ semi-ample divisors
in $d$-dimensional toric Fano varieties having at
most Gorenstein singularities.  Our main
purpose is to show that the combinatorial duality proposed by second author
agrees with the duality for Hodge numbers predicted by mirror symmetry.
It is expected that the complete verification of mirror symmetry
predictions for singular Calabi-Yau varieties $V$ in arbitrary
dimension demands considerations of
so called {\em string-theoretic Hodge numbers} $h^{p,q}_{\rm st}(V)$.
We restrict ourselves to  the string-theoretic
Hodge numbers $h^{0,q}_{\rm st}(V)$  and $h^{1,q}_{\rm st}(V)$
$(0 \leq q \leq d-r)$ which coincide with the usual Hodge numbers
$h^{0,q}(\widehat{V})$ and $h^{1,q}(\widehat{V})$ of a
$MPCP$-desingularization $\widehat{V}$ of $V$.
\end{abstract}

\thispagestyle{empty}

\newpage

\section{Introduction}

It was conjectured in \cite{bat.dual} that the polar duality for
reflexive polyhedra induces the mirror involution for Calabi-Yau
hypersurfaces in Gorenstein toric Fano varieties. The second author
has proposed a more general duality which conjecturally induces
the mirror involution also for Calabi-Yau {\em complete intersections}
in Gorenstein toric Fano varieties.
The verfication of predictions for Gromov-Witten invariants
of Calabi-Yau complete intersections in ordinary and weighted
projective spaces \cite{libgober,klemm-theisen} as well as in
some toric varieties \cite{batyrev-straten} (see also
\cite{berglund-hubsch,ell-str1,ell-str2,hosono})
can be considered as partial confirmations
of this generalized mirror construction.

If two {\em smooth} $n$-dimensional
Calabi-Yau manifolds $V$ and $W$ form a mirror
pair, then their Hodge numbers must satisfy the relation
\begin{equation}
 h^{p,q}(V) = h^{n-p,q}(W),
 \label{h.dual}
 \end{equation}
for all $0 \leq p,q \leq n$.

However, the combinatorial involutions which were constructed in
\cite{bat.dual,borisov} relate families of {\em singular}
Calabi-Yau varieties.
If $V$ is a Calabi-Yau complete intersection of semi-ample divisors in
a Gorenstein toric Fano variety ${\bf P}$, then there exists always
a partial desingularization $\pi \,: \, \widehat{V}
\rightarrow
V$ ($MPCP$-desingularization of $V$) such that:

\begin{itemize}

\item $\widehat{V}$ is again a Calabi-Yau complete intersection
of semi-ample divisors in a projective toric variety $\widehat{\bf P}$;

\item $\widehat{V}$ and $\widehat{\bf P}$ have only Gorenstein terminal
abelian quotient singularities.

\end{itemize}

If $\widehat{V}$ is smooth (this is always the case for $n \leq 3$),
we can use $\widehat{V}$ instead of
$V$ for verification of the duality (\ref{h.dual}).
In general,
we have to change the cohomology theory and to consider the so called
{\em string-theoretic Hodge numbers} $h_{\rm st}^{p,q}(V)$ for
singular $V$. For mirror pair of {\em singular} $n$-dimensional
Calabi-Yau varieties $V$ and $W$, we must have the duality for
the string-theoretic Hodge numbers:
\begin{equation}
 h^{p,q}_{\rm st} (V) = h^{n-p,q}_{\rm st}(W),
 \label{h.dual1}
 \end{equation}
for all $0 \leq p,q \leq n$.

The main properties of the string-theoretic
Hodge numbers were considered in \cite{batyrev.dais}. These numbers satisfy
the Poincar\'e duality and we have
$h^{p,q}_{\rm st}(V) = h^{p,q}_{\rm st}(\widehat{V})$ for all
$0 \leq p,q \leq n$. Moreover, one has  the following properties:

\begin{prop}
Let $h^{p,q}(\widehat{V})$
denote the usual $(p,q)$-Hodge number of $\widehat{V}$. Then

{\rm (i)}   $h^{p,q}_{\rm st}(V) = h^{p,q}(\widehat{V})$ for all
$p =0,1$ and $0 \leq q \leq n$;

{\rm (ii)}   $h^{p,q}_{\rm st}(V) = h^{p,q}(\widehat{V})$ for all
$0 \leq p,q \leq n$ if $\widehat{V}$ is smooth.
\label{string}
\end{prop}

The main purpose of this paper is to verify the duality (\ref{h.dual1})
for string-theoretic $(0,q)$ and $(1,q)$-Hodge numbers of Calabi-Yau
complete intersections $V$ and its mirror partner $W$ predicted by the
construction in \cite{borisov}. According to  \ref{string}(i),
it is sufficient to check the analgous duality for the usual Hodge numbers
of the corresponding $MPCP$-desingularizations
$\widehat{V}$ and $\widehat{W}$.
\bigskip

In section 2 we remind necessary facts from the theory of toric varieties.
Section 3 is devoted to basic properties of complete intersections in toric
varieties. In section 4 we explain the relation between Calabi-Yau complete
intersection in Gorenstein toric Fano varieties and nef-partitions of
reflexive polyhedra. In section 5 we prove the duality (\ref{h.dual1})
for $(0,q)$-Hodge numbers and give explicit formulas for them.
It turned out
to be not so easy to derive general formulas for $(1,q)$-Hodge numbers and
to prove the duality (\ref{h.dual1}) for $(1,q)$-Hodge numbers in full
generality. In sections 6 and 7 we give explicit formulas and prove the
duality only for the alternative sum of $(1,q)$-Hodge numbers; i.e.,
for the Euler characteristics of the sheaves of $1$-forms on
$\widehat{V}$ and $\widehat{W}$.
In section 8, we derive explicit formulas for
$(1,q)$-Hodge numbers of Calabi-Yau complete intersections
of {\em ample} divisors. Finally, in section 9 we establish the duality
(\ref{h.dual1}) for all  $(1,q)$-Hodge numbers of Calabi-Yau complete
intersections in projective spaces and their mirror partners.

\section{Basic notations and statements}

Let $M$ and $N = {\rm Hom}(M, {\bf Z})$ be dual free abelian groups
of rank $d$, $M_{\bf R}$ and $N_{\bf R}$ the real scalar extensions
of $M$ and $N$, $\langle *, * \rangle\; : \; M_{\bf R}
\times N_{\bf R} \rightarrow
{\bf R}$ the canonical pairing. We consider $M$ (resp. $N$)  as the
maximal lattice in $M_{\bf R}$ (resp. in $N_{\bf R}$).  We denote by
${\bf T}$ the {\em affine algebraic torus} over ${\bf C}$:
\[ {\bf T} : = {\rm Spec}\, {\bf C} \lbrack M \rbrack \cong
{\bf C} \lbrack X_1^{\pm 1}, \ldots, X_d^{\pm 1} \rbrack. \]

By a {\em lattice polyhedron} in $M_{\bf R}$ (resp. in $N_{\bf R}$)
we always mean a convex polyhedron of dimension $\leq d$
whose vertices belong to
$M$ (resp. in $N$).  The {\em relative interior} of $\Delta$ is
the  set of interior points of $\Delta$ which is considered as
a subset of the minimal ${\bf R}$-linear affine subspace
containing $\Delta$.
For any lattice polyhedron $\Delta$, we denote by $l^*(\Delta)$
the number of lattice points in the relative interior of $\Delta$.
We set  $b(\Delta) = (-1)^{{\rm dim}\, \Delta}l^*(\Delta)$ and denote by
$l(\Delta)$ the number of lattice points in $\Delta$.

A lattice polyhedron $\Delta$ defines the {\em projective toric variety} over
${\bf C}$:
\[ {\bf P}_{\Delta} = {\rm Proj}\, S_{\Delta}, \]
where $S_{\Delta}$ is the monomial subalgebra in the polynomial ring
\[ {\bf C} \lbrack X_0, X_1^{\pm 1}, \ldots, X_d^{\pm 1} \rbrack \]
spanned as ${\bf C}$-linear space
by monomials $X_0^k X_1^{m_1} \cdots X_d^{m_d}$ such that the corresponding
lattice
point $(m_1, \ldots, m_d) \in M$ belongs to $k\Delta$.
We remark that the dimension of ${\bf P}_{\Delta}$ equals
${\rm dim}\, \Delta$. If ${\rm dim}\, \Delta = d$, then ${\bf P}_{\Delta}$
can be considered as a projective compactification of ${\bf T}$. In the latter
case, the irreducible components ${\bf D}_1, \ldots, {\bf D}_n$
of ${\bf P}_{\Delta} \setminus {\bf T}$
one-to-one correspond to $(d-1)$-dimensional
 faces $\Theta_1, \ldots, \Theta_n$ of $\Delta$.
We denote by ${\bf e}_1, \ldots, {\bf e}_n$ the primitive lattice points
in $N$ which define the linear equations for the affine hyperplanes
containing $\Theta_1, \ldots, \Theta_n$ (in other words,
${\bf e}_1, \ldots, {\bf e}_n$ are primitive integral interior normal
vectors to faces $\Theta_1, \ldots, \Theta_n$).

The is another well-known definition of toric varieties ${\bf P}_{\Delta}$
via the normal fan $\Sigma = \{ \sigma_B \}$ consisting
of all rational polyhedral cones
\[ \sigma_B = {\bf R}_{\geq 0} {\bf e}_{i_1} + \cdots +
{\bf R}_{\geq 0} {\bf e}_{i_s} \subset N_{\bf R} \]
corresponding to those subsets $B = \{ i_1, \ldots, i_s \} \subset
\{1, \ldots, n\}$
for which the intersection $\Theta_{i_1} \cap \cdots \cap
\Theta_{i_s}$ is not empty. In this situation,  we also  use the
notation ${\bf P}_{\Sigma}$ for toric varieties associated with $\Sigma$.

It is known that every invertible sheaf ${\cal L}$ on any  toric variety
${\bf P}$  admits a ${\bf T}$-linearization. By this reason,
we shall consider in
this paper only ${\bf T}$-linearized invertible sheaves on toric varieties.
A ${\bf T}$-linearization of ${\cal L}$ induces the $M$-grading
of the cohomology spaces
\[ H^i({\bf P}, {\cal L}) = \bigoplus_{m \in M}
H^i({\bf P}, {\cal L})(m). \]
The convex hull $\Delta({\cal L})$ of all lattice points $m \in M$ for which
$H^0({\bf P}, {\cal L})(m) \neq 0$ will be called the {\em supporting
polyhedron
for global sections of} ${\cal L}$.

Recall the following well-known statement \cite{danilov,oda}:

\begin{theo}
There is one-to-one correspondence
\[ {\cal L} \cong {\cal O}_{\bf P}(a_1 {\bf D}_1 + \cdots + a_n {\bf D}_n)
 \leftrightarrow \varphi,\;\;
 \; a_i = \varphi({\bf e}_i), \; i = 1, \ldots, k,   \]
between ${\bf T}$-linearized invertible sheaves ${\cal L}$ on
${\bf P} = {\bf P}_{\Sigma}$ and  continious functions
$\varphi \,: \, N_{\bf R} \rightarrow {\bf R}$ which are integral
$($i.e., $\varphi(N) \subset {\bf Z})$ and ${\bf R}$-linear on
every cone $\sigma \in \Sigma$. Moreover, the supporting polyhedron
for the global sections of ${\cal L} = {\cal L}(\varphi)$ equals
\[ \Delta({\cal L}) = \{ x \in M_{\bf R} | \langle x, y \rangle \geq -
\varphi(y)\;\; \mbox{\rm for all $y \in N_{\bf R}$} \}. \]
\label{1-1}
\end{theo}

By {\em semi-ample invertible sheaf} ${\cal L}$ on a projective toric  variety
${\bf P}_{\Delta}$ we always mean an invertible
 sheaf ${\cal L}$ generated by global
sections. We  have  the following \cite{danilov,oda}:

\begin{prop}
${\cal L}$ is semi-ample if and only if the corresponding $\Sigma$-piecewise
linear $\varphi$ is upper convex.
Moreover, any  semi-ample invertible
sheaf ${\cal L}$ $($together with a ${\bf T}$-linearization$)$
on a toric variety
${\bf P}$ is uniquely determined by its  supporting polyhedron
$\Delta({\cal L})$:
\[ {\cal L} \cong {\cal O}_{\bf P}(a_1 {\bf D}_1 + \cdots + a_n {\bf D}_n),
\;\; \mbox{\rm where}\;\; a_i = - \min_{x \in \Delta({\cal L})}
\langle  x, {\bf e}_i \rangle. \]
\label{semi-ample}
\end{prop}

\begin{dfn} {\rm Let $\Delta$ and $\Delta'$ be two lattice polyhedra. Then we
call a polyhedron $\Delta'$
a {\em Minkowski summand} of $\Delta$ if there exist
a positive integer $\mu$ and  a lattice polyhedron $\Delta''$ such that
$\mu \Delta = \Delta' + \Delta''$. }
\end{dfn}

Using \ref{semi-ample1}, one easily obtains:

\begin{prop}
A lattice polyhedron $\Delta'$ is the supporting polyhedron
for global sections of a
${\bf T}$-linearized semi-ample invertible sheaf on ${\bf P}_{\Delta}$ if
and only if $\Delta'$ is a Minkowski summand of $\Delta$.
\label{semi-ample1}
\end{prop}

In the sequel, we shall use many times the following statement:

\begin{theo}
Let $D$ be a nef-Cartier divisor $($or, equivalently,
${\cal O}_{\bf P}(D)$ is a semi-ample invertible sheaf$)$
on a projective toric variety
${\bf P} = {\bf P}_{\Delta}$
$\Delta'$  the lattice polyhedron supporting the global sections
${\cal O}_{\bf P}(D)$. Then
\[ H^i({\bf P}, {\cal O}_{\bf P}(-D)) = 0, \mbox{ if $i \neq
{\rm dim}\, \Delta'$ } \]
\[ H^i({\bf P}, {\cal O}_{\bf P}(-D)) = l^*(\Delta'),
\mbox{ if $i = {\rm dim}\, \Delta'$ }. \]
In particular the Euler characteristic $\chi({\cal O}_{\bf P}(-D))$
equals $b(\Delta')$.
\label{cohomology}
\end{theo}

\noindent
{\em Proof.} Let $k  = {\rm dim}\, \Delta$. Then the invertible sheaf
${\cal O}_{\bf P}(D)$ defines the canonical morphism
\[ \pi_D \; :\;
{\bf P} \rightarrow {\bf V} =
{\rm Proj} \oplus_{n \geq 0} H^0({\bf P},
{\cal O}_{\bf P}(nD)),  \]
where ${\bf V} \cong {\bf P}_{\Delta}$ is a $k$-dimensional projective
toric variety, and
${\cal O}_{\bf P}(D) \cong \pi^*_D{\cal O}_{\bf V}(1)$. Therefore,
\[  H^i({\bf P}, {\cal O}_{\bf P}(-D)) \cong
H^i({\bf V}, {\cal O}_{\bf V}(-1)). \]
Hence $H^i({\bf P}, {\cal O}_{\bf P}(-D)) = 0$ for $i > k$.
On the other hand, $H^i({\bf V},{\cal O}_{\bf V}(-1)) = 0$ for $i < k$, and
$H^k({\bf V},{\cal O}_{\bf V}(-1)) = l^*(\Delta)$ (see \cite{danilov,khov77}).
\hfill $\Box$
\bigskip

A complex  $d$-dimensional
algebraic variety $W$ is said to have only {\em toroidal singularities}
if the $m$-adic completion of the local ring $(R, m)$ corresponding to
any  point $p \in W$ is isomorphic to the $m_{\sigma}$-completion
of a semi-group ring $(S_{\sigma}, m_{\sigma})$, where
$S_{\sigma} = {\bf C} \lbrack \sigma \cap M \rbrack$  for some $d$-dimensional
rational convex polyhedral cone $\sigma \subset M_{\bf R}$ with vertex at
$0 \in M$, and the maximal ideal $m_{\sigma} \subset S_{\sigma}$ is
generated by all non-constant monomials.
We formulate without proof the following technical statement which
is  a generalization of the classical Bertini's theorem:

\begin{theo}
Let $W$ be a complex projective algebraic variety with only toroidal
singularities, ${\cal L}$ a semi-ample invertible sheaf on $W$, $D \subset W$
the set of zeros of a generic global section of ${\cal L}$.  Then $D$ again
has only toroidal singularities. In particular, $D$ is irreducible if
${\rm dim}\, D > 0$.
\label{bertini}
\end{theo}

\section{Complete intersections}

Let $\Delta_1, \ldots, \Delta_r$ be lattice polyhedra in $M_{\bf R}$.
In this section ${\bf P}$ denotes the $d$-dimensional
toric variety ${\bf P}_{\Delta}$,  where $\Delta =
\Delta_1 + \cdots + \Delta_r$ (without loss of generality,
we assume ${\rm dim}\, \Delta = {\rm dim}\, M_{\bf R} = d$).
By \ref{semi-ample1}, the
lattice polyhedron $\Delta_i$ is the support polyhedron for global
sections of some  semi-ample invertible sheaf ${\cal L}_i$ on ${\bf P}$
$(i =1, \ldots, r
)$. We identify $H^0({\bf P}, {\cal L}_i)$
with the space of all Laurent polynomials  $f_i  \in
{\bf C}\lbrack M \rbrack$ having $\Delta_i$ as the Newton polyhedron.

\begin{dfn}
{\rm Lattice polyhedra $\Delta_1, \ldots, \Delta_r$
are called $k$-dependent if there exist $n > 0$  and $n$-element subset
\[ \{ \Delta_{i_1}, \ldots,   \Delta_{i_{n}} \} \subset
\{ \Delta_1, \ldots, \Delta_r \}, \]
such that
\[ {\rm dim}\, (\Delta_{i_1} + \cdots + \Delta_{i_{n}}) < n+k-1. \]
Lattice polyhedra which are not $k$-dependent will be called
$k$-independent.}
\end{dfn}

\begin{rem}
{\rm It follows immediately from definition that the $k$-independence of
lattice polyhedra $\Delta_1, \ldots, \Delta_r$ implies the $l$-independence
for $1 \leq  l \leq k$. }
\end{rem}

\begin{theo}
Let $Z_f \subset {\bf T}$ be a complete intersection of  $r$ affine
hypersurfaces $Z_{f_1}, \ldots, Z_{f_r}$ defined by a general system of
the equations $f_1 = \cdots = f_r = 0$ where
$f_i$ is a general Laurent polynomial with the Newton polyhedra
$\Delta_i$ $(i =1, \ldots,r)$.  Denote by $Z_i$ the closure
of $Z_{f_i}$ in ${\bf P}$ $(i =1, \ldots, r)$. Let
\[ V = Z_1 \cap \cdots \cap Z_r. \]
Then the following statements hold:

{\rm (i)} $V$ is non-empty if and only if $\Delta_1, \ldots,
\Delta_r$ are $1$-independent;

{\rm (ii)} if  $\Delta_1, \ldots, \Delta_r$ are $2$-independent,
then  $V$ is irreducible;

{\rm (iii)} if $\Delta_1, \ldots, \Delta_r$ are $k$-independent
$(k \geq 3)$, then $h^1({\cal O}_{V}) =
\cdots = h^{k-2}({\cal O}_{V}) = 0$.
\label{c}
\end{theo}

\noindent {\em Proof. }
Denote by ${\cal K}^*$ the Koszul complex
\[ {\cal O}_{\bf P}(-Z_1 - \cdots - Z_r) \rightarrow \cdots
\rightarrow \sum_{i < j} {\cal O}_{\bf P}(-Z_i -Z_j) \rightarrow
\sum_i {\cal O}_{\bf P}(-Z_i) \rightarrow {\cal O}_{\bf P}.  \]

There are two spectral sequences $'E$ and $''E$ converging to the
hypercohomology ${\bf H}^*({\bf P}, {\cal K}^*)$ (cf. \cite{griffiths}):

\[ 'E_2^{p,q} = H^p({\bf P}, {\cal H}^q({\cal K}^*)), \]
\[ ''E_2^{p,q} = H^q({\bf P}, H^p({\cal K}^*)). \]
Since ${\cal K}^*$ is which an acyclic resolution
of ${\cal O}_{V}$, we have
 ${\cal H}^q({\cal K}^*)=0$ for $q \neq r$, and
${\cal H}^r({\cal K}^*))={\cal O}_{V}$,
the first spectral sequence
degenerates and we get the isomorphisms
\[ {\bf H}^{r+p}({\bf P}, {\cal K}^*) \cong
H^{p}({\bf P}, {\cal O}_{V}) \cong
H^{p}({V}, {\cal O}_{V}
). \]
On the other hand, the second spectral sequence does not degenerate
in general. The term  ${''E}_2^{p,q}$ is the cohomology
of the bicomplex
\[ {''E}^{p,q}_1 = \bigoplus_{\{i_1, \ldots, i_{r-p}\}} H^q(
{\bf P}, {\cal O}_{\bf P}(-Z_{i_1} - \cdots - Z_{i_{r-p}})). \]
We can compute ${''E}^{p,q}_1$ using Proposition \ref{cohomology}.
The statements (i)-(iii) will follow from the consideration of
${''E}$.

(i) It is clear that $V$ is empty if and only if $Z_f$ is empty.
Assume that  $\Delta_1, \ldots, \Delta_r$ are
$1$-dependent. Then there exists $n \geq 1$ such that
$d' = {\rm dim}\, (\Delta_{i_1} + \cdots + \Delta_{i_{n}}) < n$.
This means that we can choose the coordinates $X_1, \ldots, X_d$ on
${\bf T}$ in such a way that the polynomials $f_{i_1}, \ldots,
f_{i_n}$ depends only on some $d'$ of them. Therefore, we obtain an
overdetermined system $f_{i_1} = \cdots = f_{i_n} = 0$; i.e.,
$V$ is empty.

Assume now that $V$ is empty, i.e.,
${\cal K}^*$ is acyclic. If some $\Delta_i$ is $0$-dimensional,
then  $\Delta_1, \ldots, \Delta_r$ are
$1$-dependent, and everything is proved. Otherwise, one has
the non-zero cycle in
${''E}^{r,0}_1 \cong {{''E}}^{r,0}_2 \cong
H^0({\bf P}, {\cal O}_{\bf P}) \cong {\bf C}$ which must
be killed by some next non-zero differential
\[ d_l\; :\; {''E}^{r-l,l-1}_{l} \rightarrow E^{r,0}_l\; (l \geq 2). \]
Therefore ${''E}^{r-l,l-1}_{l} \neq 0$ for some $2 \leq l \leq r$.
This implies that there exists an $l$-element subset
$\{ i_1, \ldots, i_l \} \subset \{ 1, \ldots, r \}$ such that
\[ H^{l-1}( {\bf P}, {\cal O}_{\bf P}(-Z_{i_1} - \cdots - Z_{i_l})) \neq 0.\]
Applying \ref{cohomology}, we see that there exists an $l$-element subset
\[ \{ \Delta_{i_1}, \ldots, \Delta_{i_l} \} \subset
\{ \Delta_1, \ldots, \Delta_r \} \]
such that
\[ {\rm dim}\, (\Delta_{i_1} + \cdots + \Delta_{i_{l}}) = l-1,  \]
i.e.,   $\Delta_1, \ldots,
\Delta_r$ are $1$-dependent.

(ii) Assume that $\Delta_1, \ldots, \Delta_r$ are $2$-independent.
By \ref{cohomology}, one has
\[ {''E}^{r,1}={''E}^{r-s,s-1}_1={''E}^{r-s,s}_1 =  0\; \mbox{ for
$1 \leq s \leq r$}. \]
Hence
\[ {''E}^{r,1}_l = {''E}^{r-s,s-1}_l = {''E}^{r-s,s}_l =  0\; \mbox{  for
$1 \leq s \leq r$, $l \geq 1$}. \]
So  ${\bf H}^r({\bf P}, {\cal K}^*)
\cong {\bf C} \cong H^0({V},
{\cal O}_{V})$. Therefore  ${V}$ is connected.
By \ref{bertini}, ${V}$ is irreducible.

(iii) Assume that $\Delta_1, \ldots, \Delta_r$ are $k$-independent
$(k \geq 3)$.
By the same arguments using \ref{cohomology}, we obtain
${\bf H}^r({\bf P}, {\cal K}^*)
\cong {\bf C}$, and ${\bf H}^{r+1}({\bf P}, {\cal K}^*)
= \cdots = {\bf H}^{r+k -2}({\bf P}, {\cal K}^*) = 0$. This implies
$h^1({\cal O}_{V}) =
\cdots = h^{k-2}({\cal O}_{V}) = 0$.
\hfill $\Box$
\medskip

\begin{rem}
{\rm Khovanski\^i announced (\cite{khov} p.41) that the statement \ref{c}(i)
was first discovered and proved by D. Bernshtein using
properties of mixed volumes in the following equivalent
form:}

{ The affine variety $Z_f$ is empty if and only if there exists an
$l$-dimensional affine subspace of $M_{\bf R}$ containing
affine translates of  some $l+1$  polyhedra from
$\{ \Delta_1, \ldots, \Delta_r \}$.}
\end{rem}

\begin{coro}
Assume that all lattice polyhedra $\Delta_1, \ldots, \Delta_r$
have positive dimension, $l^*(\Delta_1 + \cdots + \Delta_r ) = 1$,
$d = {\rm dim}\, (\Delta_1 + \cdots + \Delta_r)$,
and for any proper subset $\{ k_1, \ldots, k_s \} \subset
\{ 1, \ldots, r\}$ one has
$l^*(\Delta_{k_1} + \cdots + \Delta_{k_s}) = 0$.
 Then the following
statements hold:

{\rm (i)} ${V}$ is empty if and only if $d = r -1$;

{\rm (ii)} ${V}$ consists of $2$ distinct points
if and only if $d = r$;

{\rm (iii)} ${V}$ is a smooth irreducible curve of genus $1$
if and only if $d = r+1$.

{\rm (iv)} ${V}$
is a nonempty irreducible variety of dimension $d-r\geq 2$
having the property
$h^{1}({\cal O}_{V}) = \cdots =
h^{d-r-1}({\cal O}_{V}) = 0$ and $h^{0}({\cal O}_{V}) =
h^{d-r}({\cal O}_{V})= 1$ if and only if $d \geq r+2$.
\label{main.cor}
\end{coro}

\noindent
{\em Proof. } It follows from Proposition \ref{cohomology} and
our assumptions that
${''E}_1^{r,0}$ and ${''E}_1^{0,d}$ have dimension $1$ and
all remaining spaces ${''E}_1^{p,q}$ are zero.

(i) If $r = d+1$, then $Z_f$ is empty by dimension arguments.
On the other hand, if $Z_f$ is empty, then $''E^{p,q}_l$ becomes
acyclic for $ l \gg 0$. Note that the only nontrivial
one-dimensional spaces
${''E}_1^{r,0}$ and ${''E}_1^{0,d}$ can kill each other only via the
non-zero differential
\[ d_l\; :\; {''E}_l^{0,d} \rightarrow {''E}_l^{r,0} \]
where $r = d+1$ and $l = r$, i.e., $Z_f$ is empty if and only if
$r = d+1$.

Assume  $r > d+1$. Then
\[{\bf C} \cong {''E}_1^{0,d} \cong {''E}_{\infty}^{0,d} \cong
{\bf H}^d({\bf P}, {\cal K}^*).\]
 On the other hand,
the isomorphism
${\bf H}^{r+p}({\bf P}, {\cal K}^*) \cong
H^{p}(C, {\cal O}_{V})$
implies ${\bf H}^{i}({\bf P}, {\cal K}^*) = 0$
for $i < r$. Contradiction.

(ii) If $r = d$, then
\[ {\bf C}^2 \cong {''E}_1^{0,d} \oplus {''E}_1^{r,0} \cong
{''E}_{\infty}^{0,d} \oplus {''E}_{\infty}^{r,0} \cong
H^0({V}, {\cal O}_{V}). \]
Since ${V}$ is nonempty, one has ${\rm dim}\,
{V} =0$, i.e., ${V}$
consists of $2$ distinct points.

(iii)-(iv) For $r \leq d-1$, we have isomorphisms
\[ {\bf C} \cong {''E}_1^{r,0} \cong H^{0}({V},
{\cal O}_{V}), \]
\[ {\bf C} \cong {''E}_1^{0,d} \cong H^{d-r}({V},
{\cal O}_{V}), \]
and
\[ 0 \cong {''E}_1^{p,q} \cong H^{d-r}({V},
{\cal O}_{V})\;  \mbox{ if $p + q \neq
r,d$}. \]
This proves (iii)-(iv). \hfill $\Box$

\section{Calabi-Yau varieties and nef-partitions}

\begin{dfn}
{\rm A lattice polyhedron $\Delta \subset M_{\bf R}$
is called {\em reflexive} if
${\bf P}_{\Delta}$ has only Gorenstein singularities and
${\cal O}_{\bf P}(1)$ is isomorphic to the anticanonical sheaf which is
considered together with its natural ${\bf T}$-linearization. In this case,
we call   ${\bf P}_{\Delta}$  a {\em Gorenstein toric Fano variety}. }
\end{dfn}

\begin{rem}
{\rm Since the ${\bf T}$-linearized anticanonical sheaf on
${\bf P}_{\Delta}$ is isomorphic to ${\cal O}_{{\bf P}_{\Delta}}
({\bf D}_1 +
\cdots + {\bf D}_n)$, it follows from the above definition
that any reflexive polyhedron
$\Delta$ has $0 \in M$ as the unique interior lattice point. Moreover,
\[ \Delta = \{ x \in M_{\bf R} | \langle x, {\bf e}_i \rangle \geq -1, \;
i =1, \ldots, n \}, \]
where ${\bf e}_1, \ldots, {\bf e}_n$ are the primi\-tive integ\-ral inte\-rior
nor\-mal vec\-tors to codimension-$1$ faces
$\Theta_1, \ldots, \Theta_n$ of $\Delta$.
These properties of reflexive polyhedra were used in another their definition
\cite{bat.dual}.}
\end{rem}

\begin{theo} {\rm \cite{bat.dual}}
Let $\Delta$ be a reflexive polyhedron as above. Then the convex hull
$\Delta^* \subset N_{\bf R}$
of the lattice vectors ${\bf e}_1, \ldots, {\bf e}_n$ is again
a reflexive polyhedron. Moreover, $(\Delta^*)^* = \Delta$.
\end{theo}

\begin{dfn}
{\rm The lattice polyhedron $\Delta^*$ is called {\em  dual} to $\Delta$. }
\end{dfn}

Using the adjunction formula, \ref{semi-ample} and
\ref{semi-ample1}, we immediately obtain:

\begin{prop} Assume that a $d$-dimensional Gorenstein toric Fano variety
${\bf P}_{\Delta}$ contains a $(d-r)$-dimensional
complete intersection $V$ $(r<d)$
of $r$ semi-ample Cartier divisors $Z_1, \ldots, Z_r$ such that
the canonical $($or, equivalently,  dualizing$)$ sheaf of $V$
is trivial. Then there exist lattice polyhedra $\Delta_1, \ldots, \Delta_r$
such that
\[ \Delta = \Delta_1 + \cdots + \Delta_r. \]
\end{prop}

\begin{dfn}
{\rm Let $\Sigma \subset N_{\bf R}$ be the normal fan defining a Gorenstein
toric Fano variety ${\bf P}_{\Delta}$, $\varphi\, : \, N_{\bf R} \rightarrow
{\bf R}$ the integral upper convex $\Sigma$-piecewise
linear function corresponding
to the ${\bf T}$-linearized anticanonical sheaf (by \ref{1-1},
$\varphi({\bf e}_i) = 1$, $i = 1, \ldots, n$), $\Delta = \Delta_1 + \cdots +
\Delta_r$ a decomposition of $\Delta$ into a Minkowski sum of $r$ lattice
polyhedra $\Delta_j$ $(j =1, \ldots, r)$, $\varphi = \varphi_1 + \cdots +
\varphi_r$ the induced decomposition of $\varphi$ into the sum of
integral upper convex $\Sigma$-piecewise
linear functions $\varphi_j$ (see \ref{semi-ample}). Then
$\Pi(\Delta) = \{ \Delta_1, \ldots, \Delta_r \}$ is called a
{\em nef-partition of} $\Delta$ if $\varphi_j({\bf e}_i) \in \{ 0,1 \}$ for
$i = 1 , \ldots, n$, $j =1, \ldots, r$. }
\label{nef}
\end{dfn}

\begin{rem}
{\rm There are reflexive polyhedra $\Delta$
which admit  decompositions into Minkowski
sum of two lattice polyhedra, but do not admit any nef-partition. For example,
let $\Delta \subset {\bf R}^2$ be the convex hull of $4$ lattice points:
$(1,0)$, $(-1,0)$, $(0,1)$ and $(0,-1)$. Then $\Delta$ is a $2$-dimensional
reflexive polyhedron which admits a decomposition into Minkowski sum of two
$1$-dimensional polyhedra $\Delta_1$ and $\Delta_2$, where
$\Delta_1$ is the convex hull of $(-1,0)$ and $(0,-1)$, and $\Delta_2$
is the convex hull of $(0,0)$ and $(1,1)$. However, $\Delta$ does not admit
any nef-partition.}
\end{rem}

\begin{rem}
{\rm The notation in \ref{nef} are  a little bit different from
definitions and  notations in \cite{batyrev-borisov,borisov}, but they are
essentially  equivalent.}
\end{rem}

\begin{dfn}
{\rm Let $\Pi(\Delta) = \{ \Delta_1, \ldots, \Delta_r \}$ be a nef-partition
of a reflexive polyhedron $\Delta$. We define the lattice polyhedron
$\nabla_j \subset N_{\bf R}$ $( j=1, \ldots, r)$
as the convex hull of $0 \in N$ and
all lattice vectors ${\bf e}_i \in \{ {\bf e}_1, \ldots, {\bf e}_n \}$
such that $\varphi_j({\bf e}_i) = 1$. }
\end{dfn}

\begin{theo} {\rm \cite{borisov} }
The Minkowski sum $\nabla = \nabla_1 + \cdots + \nabla_r$ is a reflexive
polyhedron. Moreover, $\Pi(\nabla) = \{ \nabla_1, \ldots, \nabla_r \}$ is
a nef-partition of $\nabla$, and one has also
\[ \nabla^* = {\rm Conv} \{ \Delta_1, \ldots, \Delta_r \}, \]
\[ \Delta^* = {\rm Conv} \{ \nabla_1, \ldots, \nabla_r \}. \]
\end{theo}

\begin{dfn}
{\rm The nef-partition $\Pi(\nabla)$ is called the
{\em dual nef-partition}. }
\end{dfn}

\begin{exam}
{\rm Let $\Delta_j \subset M_{\bf R}^{(j)}$ is a reflexive polyhedron
$( j =1, \ldots, r)$, $\nabla_j : =
\Delta^*_j \subset N_{\bf R}^{(j)}$ the corresponding
dual reflexive polyhedron $(j =1, \ldots, r)$. We set
$M_{\bf R} := M_{\bf R}^{(1)} \oplus \cdots \oplus M_{\bf R}^{(r)}$,
$N_{\bf R} := N_{\bf R}^{(1)} \oplus \cdots \oplus N_{\bf R}^{(r)}$.
Then $\Pi(\Delta) = \{ \Delta_1, \ldots, \Delta_r \}$ is a nef-partition of the
reflexive polyhedron $\Delta = \Delta_1 + \cdots + \Delta_r \subset
M_{\bf R}$, and
$\Pi(\nabla) = \{ \nabla_1, \ldots, \nabla_r \}$ is the dual
nef-partition of the
reflexive polyhedron $\nabla = \nabla_1 + \cdots + \nabla_r
\subset N_{\bf R}$.}
\label{ex-split}
\end{exam}

By \cite{bat.dual},
a maximal projective triangulation ${\cal T}$
of $\Delta^* \subset N_{\bf R}$ defines a
$MPCP$-desingularization $\pi_{\cal T} \, : \, \widehat{\bf P}_{\Delta}
\rightarrow {\bf P}_{\Delta}$. The ${\bf T}$-invariant divisors
on $\widehat{\bf P}_{\Delta}$ one-to-one correspond to lattice points
on the boundary $\partial \Delta^*$.

\begin{dfn}
{\rm Denote by ${\cal V}(\Delta^*)$ the set $\partial \Delta^* \cap N$. If
$v$ is a lattice point in ${\cal V}(\Delta^*)$, then the corresponding toric
divisor on $\widehat{\bf P}_{\Delta}$ will be denoted by ${\bf D}(v)$. }
\end{dfn}

Since the anticanonical sheaf on $\widehat{\bf P}_{\Delta}$ is semi-ample,
$\Delta_1, \ldots, \Delta_r$ are supporting polyhedra  for global
sections of some semi-ample invertible sheaves
$\widehat{\cal L}_1, \ldots,\widehat{\cal L}_r$ on
$\widehat{\bf P}_{\Delta}$.

\begin{dfn}
{\rm Let  $\Pi(\Delta)$ be a nef-partition, $Z_f \subset {\bf T}$
a complete intersection of  $r$ general
affine hypersurfaces $Z_{f_1}, \ldots, Z_{f_r}$ defined by a
general system of the polynomial equations $f_1 = \cdots = f_r = 0$ where
$f_i$ is a general Laurent polynomial with the Newton polyhedra
$\Delta_i$ $(i =1, \ldots,r)$.  Denote by $\widehat{Z}_i$ the closure
of $Z_{f_i}$ in $\widehat{\bf P}_{\Delta}$ $(i =1, \ldots, r)$. Define
$V$ (resp. $\widehat{V}$) as the closure of $Z_f$ in ${\bf P}_{\Delta}$ (resp.
in $\widehat{\bf P}_{\Delta}$; i.e., $\widehat{V} =
\widehat{Z}_1 \cap \cdots \cap \widehat{Z}_r$).
If $\Pi(\nabla)$ is the dual nef-partition, then the corresponding
general complete intersection in ${\bf P}_{\nabla}$ (resp. in
$\widehat{\bf P}_{\nabla}$)  will be denoted by
$W$ (resp. $\widehat{W}$). }
\label{notation}
\end{dfn}

By the adjunction formula and \ref{bertini}, one immediately obtains:

\begin{prop}
In the above notations, assume that $V$ is nonempty and irreducible. Then
$\widehat{V}$ is an irreducible $(d-r)$-dimensional
projective algebraic variety with at most
Gorenstein terminal abelian quotient singularities and trivial canonical
class. In particular, $\widehat{V}$ is smooth if $d-r \leq 3$.
\end{prop}

In the next section we prove the following:

\begin{theo}
Let $\Pi(\Delta)$ be a nef-partition, $\Pi(\nabla)$ the dual
nef-partition, $V$ $($resp.  $W)$ corresponding to $\Pi(\Delta)$
$($resp. to $\Pi(\nabla))$ general complete intersections in
${\bf P}_{\Delta}$ $($ resp. in  ${\bf P}_{\nabla})$. Then  the following
statements hold:

{\rm (i)} $V$ is nonempty if and only if $W$ is nonempty;

{\rm (ii)} $V$ is irreducible if and only if $W$ is irreducible;

{\rm (iii)} $h^{i}({\cal O}_V) = h^{i}({\cal O}_W)$ for $0 \leq i \leq d-r$.
\label{0-hodge}
\end{theo}

\begin{coro}
Let $\Pi(\Delta)$ be a nef-partition, $\Pi(\nabla)$ the dual
nef-partition, $\widehat{V}$ $($resp.  $\widehat{W})$ corresponding
to $\Pi(\Delta)$
$($resp. $\Pi(\nabla))$ general complete intersections in
$\widehat{\bf P}_{\Delta}$ $($resp. in  $\widehat{\bf P}_{\nabla})$. Then
$\widehat{V}$ is an irreducible $(d-r)$-dimensional
projective algebraic variety with at most
Gorenstein terminal abelian quotient singularities and trivial canonical
class if and only if $\widehat{W}$ has the same properties.
\end{coro}
\medskip

\section{Semi-simplicity principle for nef-partitions}
\noindent

Let $\Delta$ be a $d$-dimensional reflexive polyhedron in $M_{\bf R}$,
$\Pi(\Delta) = \{ \Delta_1, \ldots, \Delta_r \}$ a nef-partition
of $\Delta$.

\begin{dfn}
{\rm We say that $\Pi(\Delta)$ {\em splits into a direct sum}
\[ \Pi(\Delta) = \Pi(\Delta^{(1)}) \oplus \cdots \oplus \Pi(\Delta^{(k)}) \]
if there exist convex lattice
polyhedra $\Delta^{(1)}, \ldots, \Delta^{(k)} \subset \Delta$ satisfying
the conditions

(i) $ d = {\rm dim}\,\Delta^{(1)} +  \cdots + {\rm dim}\,
\Delta^{(k)} $
and
\[ \Delta = \Delta^{(1)} + \cdots + \Delta^{(k)};\]

(ii) for $1 \leq i \leq k$, the lattice point
$0$ is cointained in the relative
interior of $\Delta^{(i)}$, $\Delta^{(i)}$ is reflexive, and
\[ \Pi(\Delta^{(i)}) = \{ \Delta_j \subset \Delta \mid
\Delta_j \subset \Delta^{(i)} \} \]
is a nef-partition of $\Delta^{(i)}$.}
\end{dfn}

\begin{dfn}
{\rm Assume that $\Pi(\Delta)$ {splits into a direct sum}
\[ \Pi(\Delta) = \Pi(\Delta^{(1)}) \oplus \cdots \oplus \Pi(\Delta^{(k)}). \]
Denote by $M^{(i)}_{\bf R}$ the minimal ${\bf R}$-linear subspace of
$M_{\bf R}$ containing $\Delta^{(i)}$ $(i =1, \ldots, k)$. We set  also
$M^{(i)} = M \cap M^{(i)}_{\bf R}$ $(i =1, \ldots, k)$. We say that
$\Pi(\Delta)$ {splits into the direct sum} {\em over {\bf Z}} if
\[ M = M^{(1)} \oplus \cdots \oplus M^{(k)}. \]}
\end{dfn}

It is easy to see the following:

\begin{prop}
Assume that a nef-partition $\Pi(\Delta)$ splits over ${\bf Z}$
into a direct sum $\Pi(\Delta^{(1)}) \oplus \cdots \oplus \Pi(\Delta^{(k)})$.
Then the Gorenstein toric Fano variety ${\bf P}_{\Delta}$ is the
product of the Gorenstein toric Fano varieties ${\bf P}_{\Delta^{(i)}}$
$(i = 1, \ldots, k)$. Moreover,
\[ V = V^{(1)} \times \cdots \times V^{(k)}, \]
\[ \widehat{V} = \widehat{V}^{(1)} \times \cdots \times \widehat{V}^{(k)}, \]
where $V^{(i)}$ $($resp. $\widehat{V}^{(i)})$ is the Calabi-Yau
complete intersections defined by the nef-partition $\Pi(\Delta^{(i)})$ in
${\bf P}_{\Delta^{(i)}}$ $($resp. in $\widehat{\bf P}_{\Delta^{(i)}})$,
$i =1, \ldots, k$.
\label{product}
\end{prop}

\begin{rem}
{\rm In \ref{ex-split} we gave a simplest  example
of a nef-partition $\Pi(\Delta)$
which splits into a direct sum $\{\Delta_1\}
\oplus \cdots \oplus \{\Delta_r\}$ over ${\bf Z}$. In general
situation, if $\Pi(\Delta)$ {splits into a direct sum}
$\Pi(\Delta) = \Pi(\Delta^{(1)}) \oplus \cdots \oplus \Pi(\Delta^{(k)})$,
then $M^{(1)} \oplus \cdots \oplus M^{(k)}$ is a sublattice of finite
index in $M$, but not necessarily the lattice $M$ itself.}
\end{rem}

\begin{exam}
{\rm Let $\Delta = \Delta_1 + \Delta_2 \subset {\bf R}^4$ where
\[ \Delta_1 = {\rm Conv}\{ (1,0,0,0), (0,1,0,0),(-1,0,0,0), (0,-1,0,0)\}; \]
\[ \Delta_2 = {\rm Conv}\{ (0,0,1,0), (0,0,0,1),(0,0,-1,0), (0,0,0,-1)\}.\]
We define the lattice $M \subset {\bf R}^4$ as
\[ M = (\frac{1}{2},\frac{1}{2},\frac{1}{2},\frac{1}{2}) + {\bf Z}^4. \]
Then $\Delta$ splits into direct sum $\Pi(\Delta) = \{\Delta_1\} \oplus
\{ \Delta_2 \}$. But it is not a splitting over ${\bf Z}$, because
$M_1 \oplus M_2$ is a sublattice of index $2$ in $M$.}
\end{exam}

\begin{dfn}
{\rm We call a nef-partition $\Pi(\Delta)$ {\em reducible} if there
exists a  subset
\[\{ k_1, \ldots, k_s \} \subset
\{1, \ldots, r \}\; ( 0< s < r) \]
such that
$\Delta_{k_1} + \ldots + \Delta_{k_s}$ contains $0$ in its relative interior.
Nef-partitions which are not reducible are called {\em irreducible}.
 }
\end{dfn}

\begin{rem}
{\rm (i) Notice that $\Delta_{k_1} + \ldots + \Delta_{k_s}$ contains $0$
in its relative interior if and only if
\[ l^*(\Delta_{k_1} + \ldots + \Delta_{k_s}) = 1. \]

(ii) It is  clear that an irreducible nef-partition has no nontrivial
splitting into a direct sum.}
\label{nul}
\end{rem}

\begin{theo}
Any nef-partition $\Pi(\Delta)$ has a unique decomposition into
direct  sum of irreducible nef-partitions.
\label{decomposition}
\end{theo}

\noindent
{\em Proof.} Assume that $\Pi(\Delta)$ is reducible. Choose
lattice polyhedra $\Delta_{k_1}, \ldots, \Delta_{k_s}$ such that
$\Delta' = \Delta_{k_1} + \ldots + \Delta_{k_s}$
contains $0$ in its relative interior. Denote
\[ \Lambda' = {\rm Conv}(\Delta_{k_1}, \ldots, \Delta_{k_s}). \]
It is
clear that $\Lambda'$ also contains $0$ in its relative interior and
has the same dimension as $\Delta_{k_1} + \ldots + \Delta_{k_s}$.
Denote by $M'_{\bf R}$ the minimal linear subspace in $M_{\bf R}$
containing $\Delta'$. Denote by $\psi_1, \ldots, \psi_r$ the
piecewise linear functions on $M_{\bf R}$ which determine the
dual nef-partition $\Pi(\nabla)$ (see \ref{nef}). We set
$\psi' = \psi_{k_1} + \cdots + \psi_{k_s}$,
$\psi'' = \psi_1 + \cdots + \psi_r - \psi' =
\psi_{j_1} + \cdots + \psi_{j_{r-s}}$. Then
\[ \Lambda' = \{ x \in M_{\bf R}' \mid \psi'(x) \leq 1 \}. \]
Since ${\psi}'$ is a non-negative integral convex
piecewise linear function having value $1$ at the relative boundary
$\partial \Lambda' \subset M_{\bf R}'$, the
polyhedron $\Lambda'$ is reflexive.

Since $\psi'$ is convex and non-negative
\[ C' = \{ x \in M_{\bf R} \mid \psi'(x) = 0 \} \]
is a convex subset of $M_{\bf R}$ containing all non-negative
linear combinations of vertices of $\Delta_{j_1}, \ldots, \Delta_{j_{r-s}}$.
One has $M_{\bf R}' \cap C' = 0$.
Assume that $0$ is not contained in the relative interior of $C'$. Then,
by separateness theorem for convex sets, there exists a non-zero element
$y' \in N_{\bf R}$ such that $\langle M_{\bf R}', y' \rangle = 0$, and
$\langle C' , y' \rangle \geq 0$. Therefore $\langle v, y' \rangle \geq 0$
for all vertices $v$ of $\Delta_1, \ldots, \Delta_r$. Contradiction.
Thus $0$ is  contained in the relative interior of $C'$. Since
${\bf R}_{\geq 0} C' \subset C'$, $C'$ is a linear subspace and
$M_{\bf R} = M_{\bf R}' \oplus C'$. We put $M_{\bf R}'' = C'$.

Define
\[ \Lambda''  = \{ x \in M_{\bf R}'' \mid \psi''(x) \leq 1 \}. \]
Then $\Lambda''=
{\rm Conv}(\Delta_{j_1}, \ldots, \Delta_{j_{r-s}})$,  and
$\Delta'':= \Delta_{j_1} \oplus \cdots \oplus \Delta_{j_{r-s}})$
contains $0$ in its relative interior. Therefore
$\Lambda''$ is also reflexive.

Thus we obtain  $\nabla^* = {\rm Conv}(\Lambda', \Lambda'')$,
${\rm dim}\, \Lambda' + {\rm dim}\, \Lambda'' = d$, and
\[ \Pi(\Delta') = \{  \Delta_{k_1}, \ldots, \Delta_{k_s} \}, \]
\[ \Pi(\Delta'') = \{ \Delta_{j_1}, \ldots, \Delta_{j_{r-s}} \} \]
are
nef-partitions  of the reflexive polyhedra $\Delta'$ and $\Delta''$
corresponding to restricitions of
$\psi_{k_1},  \ldots,  \psi_{k_s}$ (resp. of
$\psi_{j_1},  \ldots,  \psi_{j_{r-s}}$) on $M'_{\bf R}$
(resp. on $M_{\bf R}''$).
Now the statement of Theorem \ref{decomposition} follows
by induction. \hfill $\Box$

\noindent
{\bf Proof of Theorem \ref{0-hodge}.} Note that $V$ is nonempty if and
only if $h^0(V, {\cal O}_V) \neq 0$. By \ref{bertini}, $V$ is irreducible
if and only if $h^0(V, {\cal O}_V) = 1$. Therefore, (i) and (ii) follow from
(iii). By \ref{main.cor} and \ref{nul}(i),
the spectral sequence $''E^{p,q}$ which computes
the cohomology of ${\cal O}_V$ degenerates at $''E_2^{p,q}$ for all
irreducible nef-partitions. By \ref{product}, we obtain the degeneration
of $''E^{p,q}$ at $''E_2^{p,q}$ for all nef-partitions
$\Pi(\Delta)$ which split over ${\bf Z}$ into a sum of irreducible
nef-partitions.

According to Theorem \ref{decomposition}, in general situation we
have a finite covering
\[ \pi \; : \; {\bf P}_{\Delta} \rightarrow {\bf P}_{\Delta^{(1)}} \times
\cdots \times {\bf P}_{\Delta^{(k)}} \]
where the degree of $\pi$ equals the index of
$M^{(1)} \oplus \cdots \oplus M^{(k)}$ in $M$, and
$\Pi(\Delta^{(1)})$, $\ldots,$ $\Pi(\Delta^{(k)})$ are irreducible
nef-partitions.
Let
\[ \tilde{\pi}\; ;\; V \rightarrow \tilde{V} :=
V^{(1)} \times \cdots \times V^{(k)} \]
the finite covering induced by $\pi$. Then we obtain the
canonical homomorphism of the spectral sequences:
\[ \tilde{\pi}^*\; : \; ''\tilde{E}^{p,q} \rightarrow ''E^{p,q}. \]
By \ref{nul}(i), we have the canonical isomorphism
\[ ''\tilde{E}^{p,q}_1 \cong  \rightarrow ''E^{p,q}_1. \]
Since $''\tilde{E}^{p,q}$ degenerates at $''\tilde{E}^{p,q}_1$ (see the
above arguments), $''E^{p,q}$ also degenerates at $''E^{p,q}_2$.
Therefore, $h^i({\cal O}_V) = h^i({\cal O}_{\tilde{V}})$
$(i =1, \ldots, d-r)$. Analogously,
$h^i({\cal O}_W) = h^i({\cal O}_{\tilde{W}})$
$(i =1, \ldots, d-r)$. Thus, we have reduced (iii) to already known case.
\hfill $\Box$
\bigskip

\begin{coro}
Let $I = \{ 1, \ldots , r \}$. Denote by $|J|$ the cardinality of a
subset $J \subset I$. Define
\[ E(\Delta, t) = \sum_{i =0}^{d-r} h^i({\cal O}_V) t^i. \]
Then
\[ E(\Delta, t) = \sum_{J \subset I} l^*(\sum_{j \in J} \Delta_j)
t^{{\rm dim}\left( \sum_{j \in J} \Delta_j \right) - |J|}. \]
\end{coro}

By Serre duality, and using the natural  isomorphisms
$H^{i}(V, {\cal O}_V) \cong H^i(\widehat{V}, {\cal O}_{\widehat{V}})$ ($i =1,
\ldots, d-r)$,
we also obtain:

\begin{coro}
Assume that $\widehat{V}$ is nonempty and irreducible. Then
\[ h^i(\widehat{V}, {\cal O}_{\widehat{V}}) =
h^{d-r-i}(\widehat{W}, {\cal O}_{\widehat{W}}). \]
\label{o-dual}
\end{coro}

\section{$\chi(\Omega^1)$ for complete intersections}

Now we want to calculate the Euler characteristic
of $\Omega_{\widehat{{V}}}^1$ of a Calabi-Yau complete
intersection $\widehat{{V}} = \widehat{Z}_1 \cap \cdots \cap \widehat{Z}_r$
in a $MPCP$-desingularization $\widehat{\bf P} : = \widehat{\bf P}_{\Delta}$
of a Gorenstein toric Fano variety ${\bf P} : = {\bf P}_{\Delta}$.
\bigskip

One has the standard exact sequence for a complete intersection:
\[  0 \rightarrow {\cal O}_{\widehat{V}}(-\widehat{Z}_1) \oplus \cdots \oplus
 {\cal O}_{\widehat{V}}(-\widehat{Z}_r) \rightarrow
\Omega^1_{\widehat{\bf P}} \mid_{\widehat{V}}
\rightarrow \Omega^1_{\widehat{V}}
\rightarrow 0.  \]

On the other hand, for toric varieties with only quotient singularities
there exists the exact sequence \cite{danilov,oda}:

\[ 0 \rightarrow  \Omega^1_{\widehat{\bf P}}  \rightarrow
{\cal O}_{\widehat{\bf P}}^d \rightarrow \bigoplus_{v \in {\cal V}(\Delta^*)}
{\cal O}_{{\bf D}(v)} \rightarrow 0.   \]

By transversality of
$\widehat{Z}_1,  \ldots ,  \widehat{Z}_r$ to all strata on
$\widehat{\bf P}$ we can restrict the last exact
sequence on $\widehat{V}$ without obtaining additional Tor-sheaves:
\[ 0 \rightarrow  \Omega^1_{\widehat{\bf P}} \mid_{\widehat{V}} \rightarrow
{\cal O}_{\widehat{V}}^d \rightarrow \bigoplus_{v \in {\cal V}(\Delta^*)}
{\cal O}_{{\bf D}(v) \cap \widehat{V}} \rightarrow 0.   \]

Consequently, we obtain:

\begin{prop}
\[ \chi(\Omega_{\widehat{V}}^1) =  \chi({\cal O}_{\widehat{V}}^d) -
\sum_{j =1}^r \chi({\cal O}_{\widehat{V}}(-\widehat{Z}_j))  -
\sum_{v \in {\cal V}(\Delta^*)}
\chi({\cal O}_{{\bf D}(v) \cap \widehat{V}}). \]
\label{lemm1}
\end{prop}

In order to compute $\chi({\cal O}_{\widehat{V}}(-\widehat{Z}_j))$,
we consider the Koszul resolution
\[  0 \rightarrow {\cal O}_{\widehat{\bf P}}(-\widehat{Z}_1 -
\cdots - \widehat{Z}_r)
\rightarrow \cdots \]
\[ \cdots \rightarrow \sum_{j < k}
{\cal O}_{\widehat{\bf P}}( -\widehat{Z}_j -
\widehat{Z}_k ) \rightarrow
\sum_{j} {\cal O}_{\widehat{\bf P}}( -\widehat{Z}_j) \rightarrow
{\cal O}_{\widehat{\bf P}} \rightarrow {\cal O}_{\widehat{V}} \rightarrow 0.\]

Tensoring it by ${\cal O}_{\widehat{\bf P}}(-\widehat{Z}_i)$ and using
\ref{cohomology}, we get

\begin{prop}
\[ \chi({\cal O}_{\widehat{V}}(-\widehat{Z}_i))
= - \sum_j b(\Delta_i + \Delta_j) +
\sum_{j < k} b(\Delta_i + \Delta_j + \Delta_k) - \cdots \]
\[ \cdots + (-1)^r b(\Delta_i + \Delta_1 + \cdots + \Delta_r). \]
\label{lemm2}
\end{prop}

For the computation of $\chi({\cal O}_{{\bf D}(v) \cap \widehat{V}})$ we need
the following:

\begin{prop}
Let $v$ an arbitrary lattice point  in ${\cal V}(\Delta^*)$,
 $\Gamma(v)$ the minimal face  of $\Delta^*$ containing $v$. Then $\Gamma(v)$
 is a face of a polyhedron $\nabla_i$ for some $i$ $( 1 \leq i \leq r)$.
\end{prop}

\noindent
{\em Proof. } If $v$ is a vertex of $\Delta^*$, then the statement is
evident, because $\Delta^* = {\rm Conv}\{ \nabla_1, \ldots, \nabla_r \}$.

Assume that $\Gamma(v)$ is a convex hull of $k > 1$ faces of
polyhedra $\nabla_i$.
Let $v_1, \ldots, v_m$ be vertices of $\Gamma(v)$. By assumption,
$m \geq k$. Without loss of generality, we assume that
$v_1 \in \nabla_1, \ldots, v_k \in \nabla_k$.
Since $v$ belongs to relative interior of $\Gamma(v)$, there exist
positive numbers $\lambda_1, \ldots, \lambda_m$ such that
$\lambda_1 + \cdots + \lambda_m =1$ and
$\lambda_1 v_1 + \cdots + \lambda_m v_m = v$. Choose an
arbitrary vertex $w$ of the dual face $\Gamma^*(v) \subset \Delta$.
Since $\Delta = \Delta_1 + \cdots + \Delta_r$, there exist
$w_i \in \Delta_i$ such that $w = w_1 + \cdots  + w_r$.
It follows from definition  of dual nef-partitions (cf. \cite{borisov})
that
\[ \langle w_i, v_i \rangle \geq -1, \;\; 1 \leq i \leq k  \]
and
\[ \langle w_i, v_j \rangle \geq 0,\; \;  i \neq j,\; 1 \leq i \leq k, \;
1 \leq i \leq k. \]
On the other hand,
\[\langle w, v \rangle = -1 = \sum_{i=1}^m \sum_{j=1}^r  \lambda_i
\langle w_j, v_i \rangle. \]
Since $\lambda_i > 0$,
the last equality is possible only if all vertices $v_1, \ldots,
v_m$ belong to {\em the same} polyhedron $\nabla_i$ for some $i$
$(1 \leq i \leq r)$. \hfill $\Box$

\begin{coro}
\[ l(\Delta^*) = l(\nabla_1) + \cdots + l(\nabla_r) - r + 1. \]
\end{coro}
\bigskip

\begin{dfn}
{\rm Denote by $\Delta_j(v)$ the face
of $\Delta_j$ which defines the Newton polyhedron for the equation of
$\widehat{Z}_j$ in the toric
divisor ${\bf D}(v) \subset \widehat{\bf P}_{\Delta}$. }
\end{dfn}

It is easy to prove the following:

\begin{prop}
Assume that the face  $\Gamma(v) \subset \Delta^*$ is also a
face of $\nabla_i$. Let
$\Gamma^*(v)$ be the dual face of $\Delta$. Then
\[ \Delta_i(v) = \{ x \in \Delta_i \mid \langle x, v \rangle = -1 \},   \]
\[ \Delta_j(v) = \{ x \in \Delta_j \mid \langle x, v \rangle = 0\;
\mbox{ if $ j \neq i $} \},   \]
and
\[ \Delta_1(v) + \cdots + \Delta_r(v) = \Gamma^*(v). \]
\label{delta-v}
\end{prop}
\medskip

Tensoring the Koszul resolution of ${\cal O}_{\widehat{V}}$  by
${\cal O}_{{\bf D}(v)}$, we obtain:

\begin{prop}
\[ \chi({\cal O}_{{\bf D}(v) \cap C}) = - \sum_j b(\Delta_j(v)) +
\sum_{j < k} b(\Delta_j(v) + \Delta_k(v) ) - \cdots \]
\[ \cdots + (-1)^r b(\Delta_1(v) + \cdots + \Delta_r(v)). \]
\label{lemm3}
\end{prop}

Let $\nabla_i^0 = \nabla \cap {\cal V}(\Delta^*)$ $(i =1, \ldots, r)$,
and $I := \{ 1, \ldots, r \}$.
If we  rewrite
\[ \sum_{v \in {\cal V}(\Delta^*)}
\chi({\cal O}_{{\bf D}(v) \cap C}) =
\sum_{i =1}^r \sum_{v \in \nabla_i^0} \sum_{J \subset I}
(-1)^{|J| + 1} b(\sum_{j \in J} \Delta_j(v)) =   \]
\[ =  \sum_{i =1}^r \sum_{v \in \nabla_i^0}\sum_{i \not\in J \subset I}
(-1)^{|J| + 1} b(\sum_{j \in J} \Delta_j(v))
 +  \sum_{i =1}^r \sum_{v \in \nabla_i^0}\sum_{i \in J \subset I}
(-1)^{|J| + 1} b(\sum_{j \in J} \Delta_j(v)), \]
then  \ref{lemm1} and  \ref{lemm2} imply:

\begin{theo}
\[ \chi(\Omega_{\widehat{V}}^1) = \chi({\cal O}_{\widehat{V}}^d) +
\sum_{i =1}^r \sum_{J \subset I} (-1)^{|J| + 1} b(\Delta_i + \sum_{j \in J}
\Delta_j ) + \]
\[ +  \sum_{i =1}^r \sum_{v \in \nabla_i^0}\sum_{i \not\in J \subset I}
(-1)^{|J| + 1} b(\sum_{j \in J} \Delta_j(v)) + \]
\[ +  \sum_{i =1}^r \sum_{v \in \nabla_i^0}\sum_{i \in J \subset I}
(-1)^{|J| + 1} b(\sum_{j \in J} \Delta_j(v)). \]
\label{chi-formula}
\end{theo}

\section{Mirror duality for $\chi(\Omega^1)$}

We prove that the involution $\Pi(\Delta) \leftrightarrow \Pi(\nabla)$
chahge the sign of $\chi(\Omega^1)$  by $(-1)^{d-r}$ (as it would follow from
the expected
 duality $h^{1,q}(\widehat{V}) = h^{1,d-r-q}(\widehat{W})$, $0 \leq q \leq
 d-r$):

\begin{theo}
\[ \chi(\Omega_{\widehat{V}}^1) = (-1)^{d-r} \chi(\Omega_{\widehat{W}}^1). \]
\label{chi-dual}
\end{theo}

For the proof we need some preliminary statements.

\begin{prop}
Let $i \in J \subset I$,
$\Delta_i^0 = \Delta_i \cap {\cal V}(\nabla^*)$. Then a nonzero lattice point
$w$ belongs to the relative interior of the polyhedron
$\Lambda = \Delta_i + \sum_{j \in J} \Delta_j$
if and only if $w\in \Delta_i^0$ and
the zero point $0 \in N$
belongs to the relative interior of $\sum_{j \not\in J} \nabla_j(w)$.
Moreover, if this happens, then
\[ {\rm dim}\,\left(  \Delta_i + \sum_{j \in J} \Delta_j \right) +
{\rm dim}\, \left( \sum_{j \not\in J} \nabla_j(w) \right) = d. \]
\label{step1}
\end{prop}

\noindent
{\em Proof.}  First of all, let's check that the interior lattice point $w$
must belong to $\Delta_i.$

This means checking $\langle w, \nabla_k \rangle \geq 0$ for $k \neq  i$
and $\langle w, \nabla_i \rangle \geq -1$, which follows easily from the fact
that $w$ is a lattice point that belongs to $(1-\epsilon)\Lambda$ for
some small positive $\epsilon.$

The polyhedron $\Lambda - w$ is defined by the
inequalities
\[ \langle x, v \rangle \geq -2 - \langle w, v \rangle, \;
v \in \nabla_i, \]
\[ \langle x, v \rangle \geq -1 - \langle w, v \rangle,
\; v \in \nabla_j,\; j \in J,\; j \neq i,     \]
\[ \langle x, v \rangle \geq 0 - \langle w, v \rangle,
\; v \in \nabla_j,\; j \not\in J.  \]
Since  $\langle w, v \rangle \geq -1$ for  $v \in \nabla_i$, and
$\langle w, v \rangle \geq 0$ for  $v \not\in \nabla_i$,  only the
inequalities $\langle x, v \rangle \geq 0$ for $v \in \nabla_j$,
$\langle w ,
v \rangle = 0$ $(j \not\in J)$  give rise to nontrivial restrictions for the
intersection of a small neighbourhood of $0 \in M_{\bf R}$ with
$\Lambda - w$. Therefore, it remains to consider the halfspaces
defined by the inequalities $\langle x, v \rangle \geq 0$ where $ v \in
\nabla_j(w)$ $(j \not\in J)$.
Denote by $L_w$ the convex cone in $M_{\bf R}$ defined by
the inequalities
\[ \langle x, v \rangle \geq 0, \;\; \mbox{ for all
$v \in \nabla_j(w)$ and $j \not\in J$}. \]
Then  $0$ lies in the relative interior of
$\Lambda - w$ if and only if $L_w$ is a linear subspace of $M_{\bf R}$.
On the other hand, the cone
\[ C_w = \sum_{j \not\in J} {\bf R}_{\geq 0} \nabla_j(w) \]
is dual to $L_w$. Moreover, $C_w$ is a linear subspace in $N_{\bf R}$ if
and only if $0$ is contained in the relative interior of
$\sum_{j \not\in J} \nabla_j(w)$. It remains to note that a convex
cone is a linear subspaces if and only if the dual cone is a linear
subspace. In the latter case, ${\rm dim}\, L_w + {\rm dim}\, C_w = d$, i.e.,
${\rm dim}\,\left(\Delta_i + \sum_{j \in J} \Delta_j \right) +
{\rm dim}\, \left( \sum_{j \not\in J} \nabla_j(w) \right) = d$.
\hfill $\Box$

\begin{coro}
\[ \sum_{i =1}^r \sum_{J \subset I} (-1)^{|J| + 1} b(\Delta_i + \sum_{j \in J}
\Delta_j ) = (-1)^{d-r}
\left( \sum_{i =1}^r \sum_{w \in \Delta_i^0}
\sum_{i \not\in J \subset I}
(-1)^{|J| + 1} b(\sum_{j \in J} \nabla_j(w)) \right). \]
\label{term2-3}
\end{coro}

\noindent
{\em Proof.} Denote by $(-1)^{{\rm dim}\, \Theta}b^0(\Theta)$
the number of {\em nonzero}
lattice points in the relative interior of a lattice polyhedron $\Theta$.
Let $i \in J \subset I$, $J' = J \setminus \{i \}$.
Note that $0$ is in the relative interior of
$\Delta_i + \sum_{j \in J'} \Delta_j$
if and only if $0$  is in the relative interior of
$\Delta_i + \sum_{j \in J} \Delta_j$.
Since $|J| = |J'| +1$ and
\[ {\rm dim}\,\left(  \Delta_i + \sum_{j \in J'} \Delta_j \right) =
{\rm dim}\, \left( \Delta_i + \sum_{j \in J} \Delta_j \right), \]
we can have
\[\sum_{i =1}^r \sum_{J \subset I} (-1)^{|J| + 1} b(\Delta_i + \sum_{j \in J}
\Delta_j ) =
 \sum_{i =1}^r \sum_{i \in
J \subset I} (-1)^{|J| + 1} b^0(\Delta_i + \sum_{j \in J}
\Delta_j ). \]
It remains to apply Proposition \ref{step1} and the
property $|J| + |I \setminus J| = r$.
\hfill $\Box$

\begin{prop}
Let $i \in J \subset I$, and $v \in \nabla_i^0$ is a lattice point.
 Then a lattice point $w$ belongs
to the relative interior of the polyhedron
$\Lambda = \sum_{j \in J} \Delta_j(v)$
 if and only $w \in \Delta_i^0$
and the lattice point $v$
belongs to the relative interior of $\nabla_i(w) +
\sum_{j \not\in J} \nabla_j(w)$.
Moreover, if this happens, then
\[ {\rm dim}\, \left(  \sum_{j \in J} \Delta_j(v) \right) +
{\rm dim}\, \left( \nabla_i(w) + \sum_{j \not\in J} \nabla_j(w)
\right) = d -1. \]
\label{step2}
\end{prop}

\noindent
{\em Proof.}   We only need to prove the implication in one direction
and the formula for the dimensions.

 By \ref{delta-v},
$\langle  \Delta_j(v), v  \rangle = 0$ for
$j \in J$, $j \neq i$ and $\langle  \Delta_i(v), v  \rangle = -1$.
Therefore, $\langle w, v \rangle = -1$, in particular, $w \neq 0$.
Now we would like to prove that $w \in \Delta_i.$ Because of
$\Lambda \subset \Delta$, this amounts to
checking $\langle w, \nabla_k \rangle \geq 0$ for $k \neq  i.$
Suppose there exists a vertex $v'$ of $\nabla_k$ such that
$\langle w, v' \rangle = -1.$ Because $w$ lies in the relative
interior of $\Lambda$, we get $\langle \Lambda, v' \rangle = -1.$ However,
this is impossible, because $\Delta_j(v)$ contain zero if $j\neq i$,
which leads to $\Delta_i(v) \subseteq \Lambda $.

The polyhedron $\Lambda - w$ is defined by the
conditions
\[ \langle x, v' \rangle \geq -1 - \langle w, v' \rangle, \;
v' \in \nabla_j,\; j \in J,    \]
\[ \langle x, v' \rangle \geq 0 - \langle w, v' \rangle,
\; v' \in \nabla_j,\; j \not\in J,  \]
\[ \langle x, v \rangle  = -1 - \langle w, v \rangle = 0. \]
If we are interested in the neighbourhood of $0 \in M$, we
are left with the inequalities $\langle x, v' \rangle \geq 0$, where
$v'$ is either a vertex of $\nabla_j(w)$
$(j \not\in J)$, or $v'$ is the vertex of $\nabla_i(w)$, or $v'=-v$.
The zero point is in the interior
of $\Lambda - w$, hence  the convex cone $L_w$ defined by these
inequalities is a linear subspace. As a result, the dual to $L_w$ cone
\[ C_w = {\bf R}_{\geq 0}\nabla_i(w)-{\bf R}_{\geq 0}v  +
\sum_{j \not\in J} {\bf R}_{\geq 0} \nabla_j(w). \]
is also a linear subspace.

Now we use \ref{delta-v} and the equality
$\langle w, v \rangle = -1$ to conclude that
$\{ y \in C_w | \langle w, y \rangle  = 0\} $ is a linear subspace
of dimension one less which equals
\[{\bf R}_{\geq 0}\left(\nabla_i(w)-v\right)  +
\sum_{j \not\in J} {\bf R}_{\geq 0} \nabla_j(w). \]

Because $\left(\nabla_i(w)-v\right)$ and $\nabla_j(w)$ contain zero,
this implies that their sum has the same dimension and contains zero
in its interior, which proves the first statement of the proposition.

The formula for the dimensions follows from the duality of $L_w$ and $C_w.$
\hfill $\Box$

\begin{coro}
\[ \sum_{i =1}^r \sum_{v \in \nabla_i^0}
\sum_{i \in J \subset I}
(-1)^{|J| + 1} b(\sum_{j \in J} \Delta_j(v)) = \]
\[ = (-1)^{d-r} \left( \sum_{i =1}^r \sum_{w \in \Delta_i^0}
\sum_{i \in J' \subset I}
(-1)^{|J'| + 1} b(\sum_{j \in J'} \nabla_j(w)) \right). \]
\label{term4}
\end{coro}

\noindent
{\em Proof.} We set $J' = \{I \setminus J\} \cup \{i\}$.
It remains to apply Proposition \ref{step2} and the
property $|J| + |J'| = r +1$.
\hfill $\Box$
\bigskip

\noindent
{\bf Proof of Theorem \ref{chi-dual}}: It remains to
 combine the statements
 \ref{o-dual}, \ref{term2-3} and \ref{term4} with the formula in
 \ref{chi-formula}.
 \hfill $\Box$
 \bigskip

\section{Complete intersections of ample
divisors}

Our next purpose is to give explicit formulas for $(*,1)$-Hodge numbers
for $MPCP$-resolution $\widehat{V}$ of a Calabi-Yau complete intersection
of ample divisors. Notice that in this case nef-partition is
always irreducible.

First, we glue together the following two  exact sequences
\[
 0 \rightarrow {\cal O}_{\widehat{V}}(-\widehat{Z}_1) \oplus \cdots \oplus
 {\cal O}_{\widehat{V}}(-\widehat{Z}_r) \rightarrow
\Omega^1_{\widehat{\bf P}}
\mid_{\widehat{V}} \rightarrow \Omega^1_{\widehat{V}}
\rightarrow 0,  \]
\[
 0 \rightarrow  \Omega^1_{\widehat{\bf P}} \mid_{\widehat{V}} \rightarrow
{\cal O}_{\widehat{V}}^d \rightarrow \bigoplus_{v \in {\cal V}(\Delta^*)}
{\cal O}_{{\bf D}(v) \cap {\widehat{V}}} \rightarrow 0  \]
and obtain the complex
\[ {\cal Q}^*\;: \;  0 \rightarrow {\cal O}_{\widehat{V}}(-\widehat{Z}_1)
\oplus \cdots \oplus
 {\cal O}_{\widehat{V}}(-\widehat{Z}_r) \rightarrow
 {\cal O}_{\widehat{V}}^d \rightarrow \bigoplus_{v \in {\cal V}(\Delta^*)}
{\cal O}_{ {\bf D}(v) \cap \widehat{V}} \rightarrow 0 \]
whose cohomology ${\cal H}^i$ can  be  nontrivial only if $i =1$, and
\[ {\cal H}^1 ({\cal Q}^*) \cong \Omega^1_{\widehat{V}}, \]
i.e., the  hypercohomology of ${\cal Q}^*$ coincides with the cohomology of
$\Omega^1_{\widehat{V}}[1]$.

\begin{prop}
Assume that $\widehat{Z}_1, \ldots, \widehat{Z}_r$ are  nef- and big-divisors;
i.e., $\Delta_i$ is a Minkowski summand of $\Delta$ and ${\rm dim}\,
\Delta_i = {\rm dim}\, \Delta = d$ $(i = 1, \ldots, r)$. Then
\[ H^k({\cal O}_{\widehat{V}}(-\widehat{Z}_i))
= 0 \; \mbox{ for $i \neq d-r$ } \]
and
\[ h^{d-r}({\cal O}_{\widehat{V}}(-\widehat{Z}_i)) =
\sum_{J \subset I} (-1)^{r - \mid J \mid } l^*(\Delta_i + \sum_{j \in J}
\Delta_j). \]
\label{coho1}
\end{prop}

\noindent
{\em Proof. } Consider the Koszul resolution of
${\cal O}_{\widehat{V}}(-\widehat{Z}_i)$.
Then the corresponding second spectral sequence degenerates
in ${''}E_2$, because our assumptions on $\widehat{Z}_1, \ldots,
\widehat{Z}_r$
 imply  ${''}E_1^{p,q} = 0$ for $q \neq d$.
The latter immediately gives all statements. \hfill $\Box$

\begin{coro}
Let ${\widehat{V}}$ be a complete intersection such that $d-r \geq 3$.
Then, under assumptions in {\rm \ref{coho1}}, the second spectral
sequence corresponding to the complex
${\cal Q}^*$ degenerates in ${''}E_2^*$, and one obtains the
following relations
\[ h^{d-r-1}(\Omega_{\widehat{V}}^1) =
\sum_{i =1}^r h^{d-r}({\cal O}_{\widehat{V}}(-\widehat{Z}_i)) -
d - \]
\[ -  \sum_{v \in {\cal V}(\Delta^*)} \left( h^{d-r-1}(
{\cal O}_{{\bf D}(v) \cap {\widehat{V}}}) - h^{d-r-2}(
{\cal O}_{{\bf D}(v) \cap {\widehat{V}}}) \right), \]
\[ h^k(\Omega_{\widehat{V}}^1) = \sum_{v \in {\cal V}(\Delta^*)}
h^{k-2}( {\cal O}_{{\bf D}(v) \cap {\widehat{V}}})\; \; \mbox{ for
$2 \leq k \leq d-r-2$,} \]
\[ h^1 (\Omega_{\widehat{V}}^1) = \sum_{v \in {\cal V}(\Delta^*)}
h^{0}( {\cal O}_{{\bf D}(v) \cap {\widehat{V}}}) -d, \]
\[ h^0(\Omega^1_{\widehat{V}}) = h^{d-r}(\Omega^1_{\widehat{V}}) =0 .\]
\end{coro}

\begin{prop}
Assume that $\widehat{Z}_1, \ldots, \widehat{Z}_r$ are  proper
pullbacks of ample divisors on ${\bf P}_{\Delta}$; i.e.,
$\Delta_i$ and $\Delta$ are Minkowski summand of each other
$( i =1, \ldots, r)$.
Choose an element of $v \in {\cal V}(\Delta^*)$. Let $s$ be
dimension of the minimal face $\Theta \subset \Delta^*$ containing $v$.
Then the faces $\Delta_1(v), \ldots, \Delta_r(v)$ depend only
on $\Theta$ $($we denote these faces by
$\Theta^*_1, \ldots, \Theta^*_r$  $)$, and
the following statements hold:

{\rm (i)} If $d-r-s-1 >0$, then
\[ h^{d-r-s-1}({\cal O}_{{\bf D}(v) \cap {\widehat{V}}}) =
\sum_{J \subset I} (-1)^{r - \mid J \mid } l^*( \sum_{j \in J}
\Theta_j^*), \]
\[ h^0({\cal O}_{{\bf D}(v) \cap {\widehat{V}}}) =1, \]
and
\[ h^i({\cal O}_{{\bf D}(v)
\cap {\widehat{V}}}) = 0\; \; \mbox{ for all $i \neq 0, d-r-s-1$}. \]

{\rm (ii)}
If $d-r-s-1 =0$, then
\[ h^{0}({\cal O}_{{\bf D}(v) \cap {\widehat{V}}}) = 1 +
\sum_{J \subset I} (-1)^{r - \mid J \mid } l^*( \sum_{j \in J}
\Theta_j^*), \]
and
\[ h^i({\cal O}_{{\bf D}(v) \cap
{\widehat{V}}}) = 0\; \; \mbox{ for all $i \neq 0$}. \]

 {\rm (iii)}
 If  $d-r-s-1 <0$, then
${\bf D}(v) \cap {\widehat{V}}$ is empty.
\label{coho2}
 \end{prop}

\noindent
{\em Proof. } Since $\widehat{Z}_1, \ldots, \widehat{Z}_r$  are proper
pulbacks of ample divisors,
the polyhedron $\Delta^*$ is combinatorially dual to
each of $r+1$ polyhedra $\Delta$, $\Delta_1, \ldots, \Delta_r$.
By this combinatorial duality,
the faces $\Delta_i(v) \subset \Delta_i$ are dual to the face $\Theta$ and
their arbitary sums  have the same dimension $d-s-1$. In the sequel, we denote
$\Delta_i(v)$ simply by $\Theta^*_i$.

Consider the Koszul resolution of ${\cal O}_{{\bf D}(v) \cap {\widehat{V}}}$.
Then the corresponding second spectral sequence degenerates
in ${''}E_2$, because  ${''}E_1^{p,q} = 0$ for $q \neq d-s-1, 0$,
${''}E_1^{0,p} = 0$ for $p \neq r$, and ${''}E_1^{0,r} \cong {\bf C}$.
Now the statements (i)-(iii) are obvious. \hfill $\Box$

\begin{coro}
Let ${\widehat{V}}$ be a complete intersection such that $d-r \geq 3$.
 Then,
under assumptions in {\rm \ref{coho2}},
one obtains
\[ h^{d-r-1}(\Omega_{\widehat{V}}^1) =
\sum_{i =1}^r \left( \sum_{J \subset I} (-1)^{r - \mid J \mid }
l^*(\Delta_i + \sum_{j \in J}
\Delta_j)    \right) - d - \]
\[ - \sum_{\begin{array}{c} {\scriptstyle {\rm dim}\,
\Theta = 0} \\
{\scriptstyle \Theta \subset \Delta^*} \end{array}}
 \left(\sum_{J \subset I}
(-1)^{r - \mid J \mid } l^*( \sum_{j \in J}
\Theta_j^*)  \right) +
\sum_{\begin{array}{c} {\scriptstyle {\rm dim}\,
\Theta = 1} \\
{\scriptstyle \Theta \subset \Delta^*} \end{array}}
l^*(\Theta) \cdot \left(\sum_{J \subset I}
(-1)^{r - \mid J \mid } l^*( \sum_{j \in J}
\Theta_j^*)  \right),
\]
\[ h^k(\Omega_{\widehat{V}}^1) =
\sum_{\begin{array}{c} {\scriptstyle {\rm dim}\,
\Theta = d-r-k -1} \\
{\scriptstyle \Theta \subset \Delta^*} \end{array}}
l^*(\Theta) \cdot \left(\sum_{J \subset I}
(-1)^{r - \mid J \mid } l^*( \sum_{j \in J}
\Theta_j^*)  \right)\; \mbox{ for $2 \leq k \leq d-r-2$,} \]
\[ h^1 (\Omega_{\widehat{V}}^1) =
{\rm Card}\{ \mbox{ {\rm lattice points in faces of dimension} $\leq d-r-1$}
\} - d + \]
\[ + \sum_{\begin{array}{c} {\scriptstyle {\rm dim}\,
\Theta = d-r-1} \\
{\scriptstyle \Theta \subset \Delta^*} \end{array}}
l^*(\Theta) \cdot
\left(   \sum_{J \subset I} (-1)^{r - \mid J \mid } l^*( \sum_{j \in J}
\Theta_j^*) \right), \]
\[ h^0(\Omega^1_{\widehat{V}}) = h^{d-r}(\Omega^1_{\widehat{V}}) =0. \]
\label{formulas}
\end{coro}

\begin{coro}
Assume that $r =1$ and $d \geq 4$. Then the Hodge numbers
$h^{p,1}(\widehat{V})$ have the following values
\[ h^{0,1}(\widehat{V}) = h^{d-1,1}(\widehat{V}) = 0, \]
\[ h^{1,1}(\widehat{V}) = l(\Delta^*) - d -1 -
\left(\sum_{\begin{array}{c}
{\scriptstyle \Theta^* \subset \Delta^*} \\ {\scriptstyle {\rm codim}\,
\Theta^* =1} \end{array}} l^*(\Theta^*)\right) +
\left(\sum_{\begin{array}{c}
{\scriptstyle \Theta^* \subset \Delta^*} \\ {\scriptstyle {\rm codim}\,
\Theta^* =2} \end{array}} l^*(\Theta^*)\cdot l^*(\Theta) \right),   \]
\[ h^{d-2, 1}(\widehat{V}) = l(\Delta) - d -1 -
\left(\sum_{\begin{array}{c}
{\scriptstyle \Theta^* \subset \Delta^*} \\ {\scriptstyle {\rm codim}\,
\Theta^* =d}\end{array}} l^*(\Theta) \right)+
\left(\sum_{\begin{array}{c}
{\scriptstyle \Theta^* \subset \Delta^*} \\ {\scriptstyle {\rm codim}\,
\Theta^* =d-1}\end{array}} l^*(\Theta^*)\cdot l^*(\Theta) \right),  \]
 \[ h^{p,1}(\widehat{V}) =
 \sum_{\begin{array}{c}
{\scriptstyle \Theta^* \subset \Delta^*} \\ {\scriptstyle {\rm codim}\,
\Theta^* =p+1}\end{array} }l^*(\Theta^*)\cdot l^*(\Theta)
\\;\; \mbox{ for $2 \leq  p \leq d-3$ }.   \]
\end{coro}

\begin{prop}
Let ${\widehat{V}}$ be a complete intersection
of  ample divisors
$\widehat{Z}_1, \ldots, \widehat{Z}_r$ such that $d-r \geq 3$. Then
one obtains
\[ h^{d-r-1}(\Omega_{\widehat{V}}^1) =
\sum_{i =1}^r \left( \sum_{J \subset I} (-1)^{r - \mid J \mid }
l^*(\Delta_i + \sum_{j \in J}
\Delta_j)    \right) - d - \]
\[ - \sum_{\begin{array}{c} {\scriptstyle {\rm dim}\,
\Theta = 0} \\
{\scriptstyle \Theta \subset \Delta^*} \end{array}}
 \left(\sum_{J \subset I}
(-1)^{r - \mid J \mid } l^*( \sum_{j \in J}
\Theta_j^*)  \right), \]
\[ h^k(\Omega_{\widehat{V}}^1) = 0 \; \mbox{ for $k \neq 1,d-r-1$ ,} \]
\[ h^1 (\Omega_{\widehat{V}}^1) =
{\rm Card}\{ \mbox{ {\rm lattice points in faces of dimension} $\leq d-r-1$}
\} - d.  \]
\label{formulas1}
\end{prop}

\noindent
{\em Proof.} If $\widehat{Z}_1, \ldots, \widehat{Z}_r$ are ample, then
${\bf P}_{\Delta}$ has at most terminal singularities; i.e.,
$l^*(\Theta) = 0$ for all faces of $\Delta^*$ of positive dimension.
This immediately implies the formulas.
\hfill $\Box$

\begin{exam}
{\sl Complete intersection
$V_{d_1,d_2}$ of two hypersurfaces in ${\bf P}^5$.}

{\rm For two cases below, we have
\[ \sum_{ \begin{array}{c}
{\scriptstyle {\rm codim}\, \Xi =1} \\
{\scriptstyle \Xi \subset \Delta} \end{array}} l^*(\Xi) = 30. \]
\bigskip

{\bf Case 1.}
 $d_1 = d_2 = 3$. Then $l(\Delta_1) = l(\Delta_2) = 56$, $l^*(2\Delta_1) =
l^*(2\Delta_2) = 1$. Therefore,
\[ h^{2,1}(V_{3,3}) = (112 - 7) - (30 + 1+1) + 0 = 73. \]

{\bf Case 2.}
$d_1 = 2$, $d_2 = 4$. Then $l(\Delta_1) = 21$, $l(\Delta_2) = 126$,
$l^*(2\Delta_1) =0$, $l^*(2\Delta_2) = 21$. Therefore,
\[ h^{2,1}(V_{2,4}) = (147 - 7) - (30 +21+0) + 0 = 89. \]
}
\end{exam}

\begin{exam}
{\sl Complete intersection $V_{2,2,3}$ in ${\bf P}^6$.}

\[ h^{2,1}(V_{2,2,3}) = (28 + 28 + 84 - 9) - (7\cdot 6 + 1 +1 + 7 + 7 ) + 0 =
73. \]
\end{exam}

\begin{exam}
{\sl Complete intersection $V_{2,2,2,2}$ in ${\bf P}^7$.}

\[ h^{2,1}(V_{2,2,2,2}) = (4 \cdot 36 -11) - (8\cdot 7 + 3 \cdot 4 ) + 0 =
65. \]
\end{exam}

\section{Complete intersections in ${\bf P}^d$}

Consider the case when a nef-partition
$\Pi(\Delta) = \{\Delta_1, \ldots, \Delta_r\} $ defines a Calabi-Yau complete
intersection $V = V_{d_1,\ldots, d_r}$ in projective space ${\bf P}^d$.
This means that
the polyhedra $\Delta_i = d_i \Lambda $ $(i = 1, \ldots, r)$
are $d_i$-multiples of a
regular $d$-dimensional
simplex $\Lambda$ and  $d_1 + \cdots + d_r = d$.
Our purpose is to compute the $(*,1)$-Hodge numbers of the
mirror Calabi-Yau complete intersection $\widehat{W} =
\widehat{Y}_1 \cap \cdots \cap \widehat{Y}_r$  defined
by the  dual nef-partition $\Pi(\nabla) = \{ \nabla_1, \ldots, \nabla_r \}$.
In this case, $\nabla_i$ $(i = 1, \ldots, r)$ is a regular $d_i$-dimensional
simplex. If for some $i$ we have $d_i = 1$, then the Calabi-Yau complete
intersection in ${\bf P}^d$ reduces to a Calabi-Yau complete intersection in
${\bf P}^{d-1}$. So we can assume that
$d_i \geq 2$ for all $i =1, \ldots r$.

\begin{prop}
In the above situation, one has
\[ h^k({\cal O}( -
\widehat{Y}_i - \sum_{j \in J} \widehat{Y}_j)) = 0\;\; \mbox{ for all $k$ } \]
unless $J =I$, or $J \cup \{i \} = I$. Moreover
\[ h^k({\cal O}( -
\widehat{Y}_i -
\sum_{j \in I} \widehat{Y}_j)) = 0\;\; \mbox{ for $k \neq d$ },  \]
\[ h^d({\cal O}( -
\widehat{Y}_i - \sum_{j \in I} \widehat{Y}_j)) = l(\nabla_i); \]
and
\[ h^k({\cal O}(  - \sum_{j \in I} \widehat{Y}_j)) = 0\;\;
\mbox{ for $k \neq d$ },  \]
\[ h^d({\cal O}(  - \sum_{j \in I} \widehat{Y}_j)) = 1.   \]
\end{prop}

\noindent
{\em Proof.} Note that the polyhedra $\nabla_i$,  $\nabla_j$ ($ j \in J$) are
regular simplices spanning {\em linearly independent}  subspaces unless
$i \in J$, $J =I$, or $J \cup \{i \} = I$. Therefore, $\nabla_i +
 \sum_{j \in J} \nabla_j$ is a regular simplex having no lattice
 points in its relative interior
 unless $i \in J$, $J =I$, or $J \cup \{i \} = I$. If $i \in J$, then
 $\nabla_i +  \sum_{j \in J} \nabla_j$ has no lattice
 points in its relative interior since $l^*(2 \nabla_i) = 0$ for all
 $i \in I$  ($d_i = {\rm dim}\, \nabla_l \geq 2$). Hence the first statement
 follows from \ref{cohomology}.

 If $J = I$, or $j \cup \{ i \} = I$, then
 ${\rm dim}\,( \nabla_i +  \sum_{j \in J} \nabla_j)  = d$.
 By \ref{cohomology},
$h^k({\cal O}( - \widehat{Y}_i -
\sum_{j \in I} \widehat{Y}_j)) = 0$ for $k \neq d$.
The remaining statements follow from $l^*(\nabla_i + \nabla_1 +
\cdots + \nabla_r) = l(\nabla_i)$ and
$l^*(\nabla_1 + \cdots + \nabla_r) = 1$.
 \hfill $\Box$

\begin{coro}
\[ h^k({\cal O}_{\widehat{W}}(-\widehat{Y}_i)) = 0 \;
\mbox{ for $k \neq d-r$ } \]
and
\[ h^{d-r}({\cal O}_{\widehat{W}}(-\widehat{Y}_i)) = l(\nabla_i) -1. \]
\end{coro}

\noindent
{\em Proof.} Tensoring the Koszul resolution ${\cal K}^*$ of
${\cal O}_{\widehat{W}}$ by ${\cal O}(-
\widehat{Y}_i)$, we obtain the degenerated
spectral sequence $``E^{p,q}$ from which one immediately obtains
the statement.
\hfill $\Box$

\begin{prop}
For any vertex $w$ of $\nabla^*$ and any subset $J \subset I$, one has
\[ l^*(\sum_{j \in J} \nabla_j(w)) = 0, \]
unless ${\rm dim}\, \sum_{j \in J} \nabla_j(w) = 0$.
\end{prop}

\noindent
{\em Proof.}  We consider two cases:

{\bf Case 1.}  $J \neq I$.
We know that
$\sum_{j \in J} \nabla_j(w)$ is always a face of
$\sum_{j \in J} \nabla_j$. On the other hand,
the linear subspaces spanned by $\nabla_j$ are linearly independent.
Since all $\nabla_i$ are regular simplices, there is no face of
$\sum_{j \in J} \nabla_i$
of positive dimension containing a lattice point in its relative
interior.

{\bf Case 2.} $ J = I$. If there exists $ i \in I$ such that
${\rm dim}\, \nabla_i(w) = 0$, then setting $J' = J \setminus \{i \}$
we reduce all to Case 1.

The polyhedron $\nabla$ has
$d +2$ lattice points, but only those lattice points  $v \in \nabla$
which satisfy the condition $\langle w,v \rangle \in \{0,-1\}$ might appear
in $\nabla_j(w)$. Hence there exists a vertex of $\nabla$ which
does not belong to any of $\nabla_i(w)$. Therefore linear
subspaces spanned  by $\nabla_i(w)$ are linearly
independent. We again obtain the same statement,
since all $\nabla_i(w)$ are regular simplices.
\hfill $\Box$

Using the Koszul resolution, we obtain

\begin{coro}
\[ h^{k}({\cal O}_{{\bf D}(w) \cap \widehat{W}})
= 0\; \; \mbox{ for $k > 0$ } \]
\end{coro}

Therefore the second spectral sequence of the  $3$-term
complex ${\cal Q}^*$ degenerates in ${''}E_2$ and
we have

\begin{prop}
Assume that $d - r \geq 3$. Then
\[ h^k (\Omega^1_{\widehat{W}}) = 0 \; \mbox{ for $2 \leq d-r-2$ }, \]
and
\[ h^{d-r-1}(\Omega^1_{\widehat{W}}) =
- d + \sum_{i=1}^r (l(\nabla_i) -1) = 1 =
h^{1} (\Omega^1_V). \]
\end{prop}

Using the duality for the Euler characteristic from Theorem \ref{chi-dual}
\[ - h^{1} (\Omega^1_{\widehat{W}}) + (-1)^{d-r-1} h^{d-r-1}
(\Omega^1_{\widehat{W}}) = (-1)^{d-r}
\left( - h^{1} (\Omega^1_V) + (-1)^{d-r-1} h^{d-r-1}
(\Omega^1_{V}) \right), \]
we
obtain
\[ h^{1} (\Omega^1_{\widehat{W}}) = h^{d-r-1}(\Omega^1_{V}). \]
Therefore, we get the complete duality for $(1,q)$-Hodge numbers:

\begin{theo}
Let $V$ be a Calabi-Yau complete intersection of $r$ hypersurfaces
in ${\bf P}^d$ and $d -r \geq 3$, $\widehat{W}$ a $MPCP$-desingularization of
the Calabi-Yau complete intersection $W \subset {\bf P}_{\nabla}$ Then
\[ h^{q} (\Omega^1_{\widehat{W}}) = h^{d-r-q}(\Omega^1_{V})\;\;\; \mbox{ for
$0 \leq q \leq d-r$}. \]
\end{theo}


\begin{thebibliography}{99}

\bibitem{bat.dual} V. V. Batyrev,  {\em Dual polyhedra and mirror
symmetry for Calabi-Yau
hypersurfaces in toric varieties}, J. Alg. Geom., {\bf 3} (1994) 493-535.


\bibitem{batyrev-straten} V. Batyrev, D. van Straten, {\em Generalized
Hypergeometric Functions and Rational Curves on Calabi-Yau Complete
Intersections in Toric Varieties}, Preprint 1993, alg-geom 9307010,
to be published in Commun. Math. Phys.

\bibitem{batyrev-borisov} V.V. Batyrev, L.A. Borisov, {\em Dual Cones and
Mirror Symmetry for Generalized Calabi-Yau Manifolds}, to appear in
"Mirror Symmetry II" (ed. S.-T. Yau), alg-geom/9402002.

\bibitem{batyrev.dais} V. V. Batyrev, D.I. Dais, {\em Strong McKay
Corrrespondence, String-theoretic Hodge Numbers and Mirror Symmetry},
Preprint Max-Plank Inst. Bonn, MPI-94-115, alg-geom/9410001.

\bibitem{berglund-hubsch} P. Berglund, T. H\"ubsch, {\em On a Residue
Representation of Deformations, Koszul and Chiral Rings}, IASSNS-HEP-94/97,
hep-th/9411131.

\bibitem{borisov} L.A.  Borisov, {\em Towards the Mirror Symmetry for
Calabi-Yau Complete Intersections in Gorenstein Toric Fano
Varieties}, Prepint 1993, alg-geom/9310001.

\bibitem{danilov} V.I. Danilov, {\em Geometry of toric varieties},
Russ. Math. Surv., {\bf 33} (1978),97-154.

\bibitem{ell-str1} G. Ellingsrud, S.A. Str{\o}mme, {\em Counting
twisted cubic curves on general complete intersections}, alg-geom/9409006.

\bibitem{ell-str2} G. Ellingsrud, S.A. Str{\o}mme, {\em Bott's Formula and
Enumerative Geometry}, Preprint, November 1994, alg-geom/9411005.

\bibitem{griffiths} P. Griffiths, J. Harris, {\em Principles of
Algebraic Geometry}, 1978.

\bibitem{hosono} S. Hosono, A. Klemm, S. Theisen and S.T. Yau,
{\em Mirror symmetry, mirror map and application to complete intersection
Calabi-Yau spaces}, Preprint HUTMP-94-02, hep-th/9406055.

\bibitem{khov77} A. G. Khovansky, {\em Newton polyhedra and toric
varieties}, Funct. Anal. Appl. {\bf 11} (1977), 56-67.

\bibitem{khov} A. G. Khovansky, {\em Newton polyhedra and the genus of
complete intersections}, Funct. Anal. Appl. {\bf 12} (1978), 38-46.

\bibitem{klemm-theisen} A. Klemm, S. Theisen, {\em Mirror Maps and
Instanton Sums for Complete Intersections in Weighted Projective
Space}, LMU-TP 93-08, hep-th/9304034.


\bibitem{libgober} A. Libgober and J. Teitelbaum, {\em Lines on Calabi-Yau
complete intersections, mirror symmetry, and Picard-Fuchs equations},
Duke Math. J., {\bf 69}(1): 29-39, {\em Int. Math. Res. Notes} (1993).

\bibitem{oda} T. Oda, {\em Convex Bodies and Algebraic Geometry - An
Introduction to the Theory of Toric Varieties}, Ergeb. Math. Grenzgeb. (3),
vol. 15, Springer-Verlag, Berlin, Heidelberg, New York, London, Paris,
Tokyo, 1988.

\end{thebibliography}
\end{document}